\documentclass[a4paper]{article}
\usepackage[utf8]{inputenc}  
\usepackage[T1]{fontenc}
\usepackage[english]{babel}

\usepackage{lmodern} \normalfont 
\DeclareFontShape{T1}{lmr}{bx}{sc} { <-> ssub * cmr/bx/sc }{}

\usepackage{amsthm}
\usepackage{amsmath}
\usepackage{amssymb}
\usepackage{mathrsfs}
\usepackage{dsfont}
\usepackage{yhmath}
\usepackage{mathabx}

\usepackage[nottoc, notlof, notlot]{tocbibind}

\usepackage{pgf,tikz}
\usetikzlibrary{arrows}
\usetikzlibrary{patterns}

\tikzset{
xmin/.store in=\xmin, xmin/.default=-3, xmin=-3,
xmax/.store in=\xmax, xmax/.default=3, xmax=3,
ymin/.store in=\ymin, ymin/.default=-3, ymin=-3,
ymax/.store in=\ymax, ymax/.default=3, ymax=3,
}

\usepackage{subcaption}
\usepackage{geometry}
\usepackage{ifthen}

\usepackage[colorlinks=true,linkcolor=blue,citecolor=green]{hyperref}

\newcommand{\R}{\mathbb{R}}
\newcommand{\N}{\mathbb{N}}
\newcommand{\Z}{\mathbb{Z}}

\newcommand{\C}{\mathbb{C}}
\newcommand{\T}{\mathbb{T}}
\renewcommand{\O}{{\cal O}}

\def\bsys{\left\{\begin{array}}
\def\esys{\end{array}\right.}

\newcommand{\eps}{\varepsilon}

\renewcommand{\emptyset}{\varnothing}
\newcommand{\de}{\,\mathrm{d}}
\renewcommand{\i}{\mathrm{i}}
\newcommand{\intg}[2]{\int\limits_{#1}^{#2}}
\renewcommand{\sl}{\mathrm{sl}}
\newcommand{\Hom}{\mathrm{Hom}}

\newtheorem{theorem}{Theorem}[section]
\newtheorem{lemma}[theorem]{Lemma}
\newtheorem{proposition}[theorem]{Proposition}
\newtheorem{corollary}[theorem]{Corollary}
\newtheorem{definition}[theorem]{Definition}
\newtheorem{remark}[theorem]{Remark}
\newtheorem{example}[theorem]{Example}
\newtheorem{nexample}[theorem]{Non-example}

\newcommand{\lem}[2][None]
{\begin{lemma}\ifthenelse{\equal{#1}{None}}{}{\label{#1}}
#2
\end{lemma}}

\newcommand{\prop}[2][None]
{\begin{proposition}\ifthenelse{\equal{#1}{None}}{}{\label{#1}}
#2
\end{proposition}}

\newcommand{\cor}[2][None]
{\begin{corollary}\ifthenelse{\equal{#1}{None}}{}{\label{#1}}
#2
\end{corollary}}

\newcommand{\theo}[2][None]
{\begin{theorem}\ifthenelse{\equal{#1}{None}}{}{\label{#1}}
#2
\end{theorem}}

\newcommand{\theonm}[2]{\smallskip \noindent \textbf{Theorem #1} {\it #2} \smallskip}

\newcommand{\defi}[2][None]
{\begin{definition}\ifthenelse{\equal{#1}{None}}{}{\label{#1}}
#2
\end{definition}}

\newcommand{\rem}[2][None]
{\begin{remark}\ifthenelse{\equal{#1}{None}}{}{\label{#1}}
#2
\end{remark}}

\newcommand{\expl}[2][None]
{\begin{example}\ifthenelse{\equal{#1}{None}}{}{\label{#1}}
#2
\end{example}}

\newcommand{\nexpl}[2][None]
{\begin{nexample}\ifthenelse{\equal{#1}{None}}{}{\label{#1}}
#2
\end{nexample}}

\newcommand{\demo}[1]{\begin{proof} #1
\end{proof}}

\newcommand{\demode}[2]{\begin{proof}[Proof #1] #2
\end{proof}}

\newcounter{exerc}

\newcounter{quest}

\newcommand{\qst}[2][None]{\refstepcounter{quest}\ifthenelse{\equal{#1}{None}}{}{\label{#1}} \smallskip \noindent {\bf Question \thequest\;} {#2} \smallskip}

\newcommand{\exo}[2][None]{\refstepcounter{exerc}\ifthenelse{\equal{#1}{None}}{}{\label{#1}} \smallskip \noindent {\bf Exercise \theexerc\;} {#2} \smallskip}

\title{Examples and non-examples of polyhedral Kähler surfaces}
\author{Cécile Gachet}

\begin{document}

\maketitle

\begin{abstract}
Polyhedral Kähler surfaces are a class of complex surfaces, which are flat everywhere except on a two-dimensional skeleton. They are defined as a generalisation of the ``gluing a polygon side by side'' construction of flat Riemann surfaces.

In this article, we introduce two classes of polyhedral Kähler surfaces with trivial holonomy: products of polyhedral Kähler curves with zero holonomy and ramified coverings of tori, and prove that none of these classes is contained in the other. Existence of other types of polyhedral Kähler surfaces is still open.
\end{abstract}

\tableofcontents

\newpage

\section{Introduction}

\subsection{Definition and examples of $PK$ manifolds}

There is a well-known construction of flat Riemann surfaces with conical singularities, occurring by gluing two-by-two sides of same lengths of a polygon $P\subset\C$. The resulting topological surface $S$ naturally inherits a complex structure, which is flat except at a finite number of singular points: charts are trivial at neighbourhoods of points not being vertices of $P$, and any vertex $x$ of $P$ has a neighbourhood $U_x\subset S$ isometric to the cone $C_{\alpha}\subset\C$ centred at zero of conical angle $2\pi\alpha$, so that the map $z\in C_{\alpha}\simeq U_x\longmapsto z^{\alpha}\in\C$ is a chart.

This point of view on Riemann surfaces is given several motivations in the survey \cite{Zorich} by A. Zorich. One of them is the definition of the so-called \emph{Teichmüller geodesic flow} on the moduli space of genus $g$ Riemann surfaces (and even on the Teichmüller spaces, because the flow preserves a given choice of representatives of the fundamental group). It is defined as follows: fix a polygon with sides which are two-by-two identified by translations, so that you can glue it to a flat Riemann surface $X_0$. Fixing a real direction in this polygon, scaling in this direction and gluing the new polygon according to the old gluing-pattern, we get a new flat Riemann surface of the same topological type as $X_0$. This defines our Teichmüller geodesic flow. Its orbits and ergodic properties are studied in geometry as well as in dynamical systems.

The initial gluing-procedure can be generalised to higher dimensions to construct so-called \emph{polyhedral manifolds}.

\defi{Consider $M$ a $2n$-dimensional $\R$-smooth manifold together with some fixed simplicial decomposition, and endow it with a flat metric obtained by isometrical gluing of choices of flat metrics on each $2n$-dimensional simplex. Then $M$ is called a \emph{polyhedral manifold}.}

From now on, let us consider for convenience only connected manifolds.

\defi{If $M^{2n}$ is a polyhedral manifold, let $M_s$ be the union of all its simplices of codimension at least 2. Then $M\setminus M_s$ clearly has a complex structure, so that the \emph{holonomy} of $M$ is well-defined as the holonomy of $M\setminus M_s$. This holonomy is called \emph{unitary} if it is a subgroup of $U(n)$.}

For $n=1$, we were lucky enough to get an additional complex structure on $M$. This does not arise in general, for $M$ may have a non-unitary holonomy: $U(1)=SO(2)$ but $U(n)\subsetneq SO(2n)$ in higher dimension. 

Actually, a more careful study sheds light on another obstruction. Let $M^{2n}$ be a polyhedral manifold of unitary holonomy, $M_s$ its 2-skeleton. Let us introduce a few useful concepts for studying the singular locus of the flat metric, $M_s$. First, if $x\in M_s$ does not belong to any simplex of codimension greater or equal to 3 ($x$ is in a face $F_x$ of codimension 2 of the simplicial decomposition), then there is a neighbourhood $x\in U_x\subset M$ isometric to $C_{\alpha}\times\R^{2n-2}$. We call $2\pi\alpha$ the \emph{conical angle} at $x$. In particular, if $M$ is $PK^{0}$, all such conical angles are in $2\pi\Z$. 
Moreover, this isometry carries the natural parallel complex structure of $U_x\cap (M\setminus M_s)$ to a parallel complex structure on $C_{\alpha}^*\times\R^{2n-2}$, which can be extended by continuity to $C_{\alpha}\times\R^{2n-2}$. If $\{0\}\times\R^{2n-2}$ is stable under this complex structure, the face $F_x$ is said to have a \emph{complex direction}.

\defi{A polyhedral manifold $M^{2n}$ is said to be \emph{polyhedral Kähler (PK)} if it has unitary holonomy and each face of codimension 2 has a complex direction. It is called a \emph{$PK$ curve} if $n=1$, a \emph{$PK$ surface} if $n=2$. We may also encounter polyhedral Kähler manifolds with trivial holonomy, being referred to as $PK^0$ manifolds.}

\expl{Let $n\ge 2$. A covering of the torus $T^4$ ramified with conical angle $2\pi n$ over the complex curve $\{(x=a+\i b,y=c+\i d)\in T^4\mid y=0\}$ is a $PK$ manifold. It also holds if ramified over the complex curve $\{(x,y)\in T^4\mid x=y\}$. If ramified over $\{(a,b,c,d)\in T^4\mid b=d=0\}$ however, it is not polyhedral Kähler, for the ramified locus has no complex direction.}

More generally, the following result holds. Its proof is postponed to Section {\bf 3.2}, where various results about complex tori are recalled.

\prop[covtispk]{Let $\T$ be a two-dimensional complex torus and $C_1,\ldots,C_n$ complex curves with fixed complex directions in $\T$. Assume that no three of them have a common intersection point. Then the ramified cover of $\T$ over these curves is a $PK^0$ surface.}

\expl{We will see in Section {\bf 2.1} that $PK^0$ curves are exactly flat Riemann surfaces with trivial holonomy and a finite number of singular points. A product of two flat Riemann surfaces with trivial holonomy and a finite number of singular points is a $PK^0$ surface.}

\rem{There is \emph{a priori} no way to scale a given $PK^0$ surface $X$, let alone to define an analogous of the Teichmüller geodesic flow on moduli spaces of $PK^0$ surfaces. That is a motivation for us to study some families of $PK^0$ surfaces on which natural scalings may arise. For instance, $PK^0$ product of curves definitely have such scalings. On the other hand, even a $PK^0$ surface $X$ arising as a ramified covering $\pi\,:\, X\to B$ of a scalable $PK^0$ surface $B$ may have no scaling induced by $B$. A scaling of $B$ \emph{lifting up} to $X$ is indeed defined as a scaling of $B$ compatible with the $PK^0$ structure considered on $X$; for that, it should send any irreducible component of $\pi(X_s)$ to a curve in $B$ which still has a complex direction.}

\rem{If a scaling of $B$ is not compatible with one $PK^0$ structure of $X$ (in the sense of the previous remark), it may still be compatible with another (biholomorphic but not isometric) $PK^0$ structure on $X$. Unfortunately, we do not understand relations between various biholomorphic $PK^0$ structures on a given $PK^0$ manifold.}

\subsection{Complex Kähler structure and flat metric on $PK$ surfaces or curves}

We already mentioned a natural complex structure arising on any $PK^0$ curve. Actually, $PK^0$ surfaces are trickier to handle; let us just quote here a theorem of \cite{Panov}:

\theonm{(by Dima Panov)}{On any polyhedral Kähler surface $M$, a natural parallel complex structure on $M\setminus M_s$ can be extended to $M$ by charts on points of $M_s$, in order to obtain a smooth complex structure on $M$. This complex flat structure is Kähler. Moreover, $M_s$ is a union of complex curves for this structure.}

We do not need this theorem in Section {\bf 2}: $PK^0$ surfaces arising as a product of two $PK^0$ curves have indeed an explicit complex structure. It will appear to be mostly important in Section {\bf 3}, significantly enough the third part about the structure of $M_s$.

Let us denote by $\Omega^1(X)$ the sheaf of holomorphic 1-forms on a complex manifold $X$.
Combining complex structure and flat metric on $PK^0$ surfaces and curves yields to the following definition:

\defi{A holomorphic 1-form on an $n$-dimensional $PK^0$ manifold is said to be \emph{compatible with the flat structure} if at any non-singular point, it is locally pulled back from a form of constant coefficients $a_1\de z_1+\ldots+ a_n\de z_n$ on $\C^n$ by the isometric local chart.}

\rem{Particularly, on a $PK^0$ curve $X$, the set of all holomorphic 1-forms compatible with the flat structure is 1-dimensional in $H^0(\Omega^1(X))$, which is itself of dimension $g(X)$. Taking this compatibility with the flat structure into account will appear as a nice way to distinguish a special abelian differential (up to multiplicative constant). It will be used thoroughly in the proof of theorem B.} 

At last, we should say that the ``Kähler part'' of D. Panov's theorem is not proven in \cite{Panov}. For sake of self-contentedness, let us sketch here an argument for proving due to Misha Verbitsky and Dima Panov. It uses the following lemma:

\lem{On a non-Kähler surface, squares of (1,1)-forms are non-positive.}

The proof of this lemma is to be found in an upcoming book by Misha Verbitsky and Liviu Ornea about the classification of complex surfaces. 

In this article, we actually just need to prove that a $PK^0$ surface can not be a Hopf surface (Section {\bf 3.3}), hence finding a non-zero (1,1)-form is already enough for our purpose (even if it has some degenerations).

\demode{of the existence of a (1,1)-form of positive square on a $PK^0$ surface}{Let $M$ be a $PK^0$ surface with singular locus $M_s$ of real codimension at least 2 and $S$ singular locus of codimension 4. Let $\omega$ be the Kähler form associated with the flat metric on $M\setminus M_s$. 

For $p\in M_s\setminus S$, $U_p$ a chart neighbourhood of $p$ and $V$ a flat open chart set intersecting $U_p$, with flat coordinates $z_1,z_2$. There are some coordinates $(z,z_1)$ on $U_p$, $(z_3,z_1)$ on $\C^2$ in which there is a chart $(z,z_1)\in U_p \to(z^n,z_1)\in\C^2$ centered at $p=(0,0)$ with $M_s=\{z=0\}$ locally. We can also write $z_3:=az_1+bz_2=0$ for some constants $a,b\in\C$, say $b\ne 0$. Then on $U_x\cap V$, $\de z_1 = \de z_1$ and $\de z_2 = b^{-1}(nz^{n-1}\de z -a\de z_1)$. Since 
$$\omega=\frac{\i}{2}(\de z_1\wedge\de\overline{z_1}+\de z_2\wedge\de\overline{z_2})$$ on $V$, rewriting in coordinates $(z,z_1)$ on $U_x\cap V$ extends it holomorphically to $\omega= \frac{\i}{2}(1+|ab^{-1}|^2)\de z_1\wedge\de\overline{z_1}$ on $M_s\cap U_x$. This does not depend of the choice of $V$ since $M$ has trivial holonomy, so that $\de z_1$ can be chosen globally on $M\setminus M_s$.

For $p\in S$ now, we can extend $\omega$ by Hartogs theorem (in the basis of local (1,1)-forms, $\omega$ is indeed given by a finite number of holomorphic functions defined everywhere except at the (real codimension 4) point $p$).

Hence, $\omega$ extends to a global closed (1,1)-form on $M$. It has some degenerations, so it is not Kähler itself; nevertheless its square is still the local volume form almost everywhere, so it equals 2.}


\subsection{Results}

Remember that a $PK^0$ curve $X$ has an Albanese variety: $\mathrm{Alb}(X)=H^0(\Omega^1(X))^*/\iota(H_1(X,\Z))$ where $$\iota\,:\,[\gamma]\in H_1(X,\Z)\mapsto \left(\omega\mapsto\intg{\gamma}{}\omega\right)\in H^0(\Omega^1(X))^*$$ is an embedding. This $\mathrm{Alb}(X)$ is a complex torus, and fixing $x_0\in X$, we have a map $$x\in X\mapsto x-x_0\in\mathrm{Pic}^0(X)\simeq\mathrm{Alb}(X).$$

Even with the $PK$ assumption though, there is no reason why this map should be a (-n eventually ramified) covering. Actually, we prove in Section {\bf 2} that:

\theonm{A}{A $PK^0$ curve is a ramified covering of a torus if and only if there are two directions $V$ and $H$ such that $V$-periods (respectively $H$-periods) form a discrete subgroup of $\R$.}

\rem{This is actually equivalent to the fact that for any two distinct directions $V$ and $H$, $V$-periods (resp. $H$-periods) form a discrete subgroup of $\R$, or also to the fact that periods form a lattice in $\C$.}

This idea enables us to construct $PK^0$ surfaces which are products of two curves but not ramified coverings of tori, by theorem B of Section {\bf 2.3}:

\theonm{B}{A product of two $PK^0$ curves $X=C_1\times C_2$ is a ramified covering of a torus $\T$ if and only if $C_1$ and $C_2$ both are ramified coverings of tori.}

In Section {\bf 3}, we consider the converse problem: which $PK^0$ surfaces are products of two curves? The answer depends actually on the set of complex directions of $X$ realised by codimension 2 faces. After a few definitions postponed to Section {\bf 3}, we will be able to prove the following result:

\theonm{C}{Let $X$ be a $PK^0$ surface. If it has at most two relevant complex directions, it is a product of curves. Moreover, there are $PK^0$ surfaces with three relevant complex directions not being products of curves.} 

\section{Does a $PK^0$ product of curves cover a torus?}

\subsection{Characterising $PK^0$ curves}

All starts with the following correspondence.

\prop[flatRS]{A $PK^0$ curve is the same as a zero holonomy flat Riemann surface with a finite number of singular points.}

The direct inclusion has already been proved.

Let $X$ be a zero holonomy flat Riemann surface with a finite number of singular points. Since the surface has zero holonomy, we can fix a direction in the tangent space at each point of $X$; let's call it the \emph{north}. Completing this unitary vector field into an indirect orthonormal basis of $TX$, we can define an orthogonal unitary vector field \emph{east}. For all $x\in X$, we denote by $l_x\,:\,\R\to X$ the geodesic path of northern direction such that $l_x(0)=x$ and $m_x$ the analogous eastern path.

To prove that $X$ is a $PK^0$ curve, we are going to use the following result.

\lem[exclusive]{There is a compact geodesic arc of eastern direction $S$ such that for all $x\in X$, the path $l_x$ crosses $S$.}

\demo{Actually, take any $p\in X$ and the geodesic arc of eastern direction $S=m_p([-T,T])$, for some $T$ to be chosen in a while. Let $U=\{x\in X\mid l_x\mbox{ crosses }S\}$. It contains $S$. For a well-chosen $T$, $U$ even contains a neighbourhood of $S$. Indeed, $m_p(\R)\subset X$ is compact, so finitely covered by open sets $U_1,\ldots,U_N$ of $X$ trivialising the north-east grid and of eastern diameter smaller than 2. Without loss of generality, $m_p(\R)\cap U_j\ne\emptyset$ for all $j$. Since $m_p(\R)=\bigcup_{n>0}m_p([-n,n])$, there is some big enough $n$ such that $m_p([-n,n])\cap U_j\ne\emptyset$ for all $j$. Then $T=n+1$ is fine.

\begin{figure}
\centering
\definecolor{qqffqq}{rgb}{0.0,1.0,0.0}
\definecolor{xdxdff}{rgb}{0.49019607843137253,0.49019607843137253,1.0}
\definecolor{uuuuuu}{rgb}{0.26666666666666666,0.26666666666666666,0.26666666666666666}
\definecolor{qqqqff}{rgb}{0.0,0.0,1.0}
\definecolor{ffqqqq}{rgb}{1.0,0.0,0.0}
\begin{tikzpicture}[line cap=round,line join=round,>=triangle 45,x=0.6cm,y=0.6cm,scale=1.3]
\clip(-2,-2) rectangle (12,3);
\draw(0.0,0.0) circle (0.848528137423857cm);
\draw(5.0,-0.0) circle (0.848528137423857cm);
\draw(10.0,-0.0) circle (0.848528137423857cm);
\draw(0.1,2) node{$U_1$};
\draw(5.1,2) node{$U_2$};
\draw(10.1,2) node{$U_3$};
\begin{scope}
\clip(0.0,0.0) circle (0.848528137423857cm);
\draw (-1.2000000000000002,-5.320000000000002) -- (-1.2000000000000002,6.3);
\draw (-0.8000000000000002,-5.320000000000002) -- (-0.8000000000000002,6.3);
\draw (-0.40000000000000013,-5.320000000000002) -- (-0.40000000000000013,6.3);
\draw (-0.0,-5.320000000000002) -- (-0.0,6.3);
\draw (0.4,-5.320000000000002) -- (0.4,6.3);
\draw (0.8,-5.320000000000002) -- (0.8,6.3);
\draw (1.2000000000000002,-5.320000000000002) -- (1.2000000000000002,6.3);
\draw (1.6,-5.320000000000002) -- (1.6,6.3);
\draw (2.0,-5.320000000000002) -- (2.0,6.3);
\draw (2.4,-5.320000000000002) -- (2.4,6.3);
\draw (2.8,-5.320000000000002) -- (2.8,6.3);
\draw (3.1999999999999997,-5.320000000000002) -- (3.1999999999999997,6.3);
\draw (3.5999999999999996,-5.320000000000002) -- (3.5999999999999996,6.3);
\draw (4.0,-5.320000000000002) -- (4.0,6.3);
\draw (4.4,-5.320000000000002) -- (4.4,6.3);
\draw (4.800000000000001,-5.320000000000002) -- (4.800000000000001,6.3);
\draw (5.200000000000001,-5.320000000000002) -- (5.200000000000001,6.3);
\draw (5.600000000000001,-5.320000000000002) -- (5.600000000000001,6.3);
\draw (6.0,-5.320000000000002) -- (6.0,6.3);
\draw (6.4,-5.320000000000002) -- (6.4,6.3);
\draw (6.800000000000001,-5.320000000000002) -- (6.800000000000001,6.3);
\draw (7.200000000000001,-5.320000000000002) -- (7.200000000000001,6.3);
\draw (7.600000000000001,-5.320000000000002) -- (7.600000000000001,6.3);
\draw (8.0,-5.320000000000002) -- (8.0,6.3);
\draw (8.4,-5.320000000000002) -- (8.4,6.3);
\draw (8.8,-5.320000000000002) -- (8.8,6.3);
\draw (9.200000000000001,-5.320000000000002) -- (9.200000000000001,6.3);
\draw (9.600000000000001,-5.320000000000002) -- (9.600000000000001,6.3);
\draw (10.0,-5.320000000000002) -- (10.0,6.3);
\draw (10.4,-5.320000000000002) -- (10.4,6.3);
\draw (10.8,-5.320000000000002) -- (10.8,6.3);
\draw (11.200000000000001,-5.320000000000002) -- (11.200000000000001,6.3);
\draw (11.600000000000001,-5.320000000000002) -- (11.600000000000001,6.3);
\draw (12.0,-5.320000000000002) -- (12.0,6.3);
\draw [domain=-4.3:18.7] plot(\x,{(-2.0-0.0*\x)/1.0});
\draw [domain=-4.3:18.7] plot(\x,{(-1.4-0.0*\x)/1.0});
\draw [domain=-4.3:18.7] plot(\x,{(-0.7999999999999999-0.0*\x)/1.0});
\draw [domain=-4.3:18.7] plot(\x,{(-0.19999999999999996-0.0*\x)/1.0});
\draw [domain=-4.3:18.7] plot(\x,{(--0.4-0.0*\x)/1.0});
\draw [domain=-4.3:18.7] plot(\x,{(--1.0-0.0*\x)/1.0});
\draw [domain=-4.3:18.7] plot(\x,{(--1.6-0.0*\x)/1.0});
\draw [line width=1.2000000000000002pt,color=cyan,domain=-4.3:18.7] plot(\x,{(--0.2-0.0*\x)/1.0});
\draw [line width=1.2000000000000002pt,color=ffqqqq,domain=-4.3:18.7] plot(\x,{(-0.5-0.0*\x)/1.0});
\draw [line width=1.2000000000000002pt,color=ffqqqq,domain=-4.3:18.7] plot(\x,{(-0.65-0.0*\x)/1.0});
\draw [line width=1.2000000000000002pt,color=cyan] (3.677124344467705,-0.5)-- (6.322875655532295,-0.5);
\draw [line width=1.2000000000000002pt,color=cyan] (9.62,-0.5)-- (10.52,-0.5);
\draw [line width=1.2000000000000002pt,color=qqffqq] (-2.24,0.2)-- (9.48,0.2);
\draw [line width=1.2000000000000002pt,color=qqqqff] (-0.94,0.2)-- (1.04,0.2);
\end{scope}
\begin{scope}
\clip(5.0,-0.0) circle (0.848528137423857cm);
\draw (-1.2000000000000002,-5.320000000000002) -- (-1.2000000000000002,6.3);
\draw (-0.8000000000000002,-5.320000000000002) -- (-0.8000000000000002,6.3);
\draw (-0.40000000000000013,-5.320000000000002) -- (-0.40000000000000013,6.3);
\draw (-0.0,-5.320000000000002) -- (-0.0,6.3);
\draw (0.4,-5.320000000000002) -- (0.4,6.3);
\draw (0.8,-5.320000000000002) -- (0.8,6.3);
\draw (1.2000000000000002,-5.320000000000002) -- (1.2000000000000002,6.3);
\draw (1.6,-5.320000000000002) -- (1.6,6.3);
\draw (2.0,-5.320000000000002) -- (2.0,6.3);
\draw (2.4,-5.320000000000002) -- (2.4,6.3);
\draw (2.8,-5.320000000000002) -- (2.8,6.3);
\draw (3.1999999999999997,-5.320000000000002) -- (3.1999999999999997,6.3);
\draw (3.5999999999999996,-5.320000000000002) -- (3.5999999999999996,6.3);
\draw (4.0,-5.320000000000002) -- (4.0,6.3);
\draw (4.4,-5.320000000000002) -- (4.4,6.3);
\draw (4.800000000000001,-5.320000000000002) -- (4.800000000000001,6.3);
\draw (5.200000000000001,-5.320000000000002) -- (5.200000000000001,6.3);
\draw (5.600000000000001,-5.320000000000002) -- (5.600000000000001,6.3);
\draw (6.0,-5.320000000000002) -- (6.0,6.3);
\draw (6.4,-5.320000000000002) -- (6.4,6.3);
\draw (6.800000000000001,-5.320000000000002) -- (6.800000000000001,6.3);
\draw (7.200000000000001,-5.320000000000002) -- (7.200000000000001,6.3);
\draw (7.600000000000001,-5.320000000000002) -- (7.600000000000001,6.3);
\draw (8.0,-5.320000000000002) -- (8.0,6.3);
\draw (8.4,-5.320000000000002) -- (8.4,6.3);
\draw (8.8,-5.320000000000002) -- (8.8,6.3);
\draw (9.200000000000001,-5.320000000000002) -- (9.200000000000001,6.3);
\draw (9.600000000000001,-5.320000000000002) -- (9.600000000000001,6.3);
\draw (10.0,-5.320000000000002) -- (10.0,6.3);
\draw (10.4,-5.320000000000002) -- (10.4,6.3);
\draw (10.8,-5.320000000000002) -- (10.8,6.3);
\draw (11.200000000000001,-5.320000000000002) -- (11.200000000000001,6.3);
\draw (11.600000000000001,-5.320000000000002) -- (11.600000000000001,6.3);
\draw (12.0,-5.320000000000002) -- (12.0,6.3);
\draw [domain=-4.3:18.7] plot(\x,{(-2.0-0.0*\x)/1.0});
\draw [domain=-4.3:18.7] plot(\x,{(-1.4-0.0*\x)/1.0});
\draw [domain=-4.3:18.7] plot(\x,{(-0.7999999999999999-0.0*\x)/1.0});
\draw [domain=-4.3:18.7] plot(\x,{(-0.19999999999999996-0.0*\x)/1.0});
\draw [domain=-4.3:18.7] plot(\x,{(--0.4-0.0*\x)/1.0});
\draw [domain=-4.3:18.7] plot(\x,{(--1.0-0.0*\x)/1.0});
\draw [domain=-4.3:18.7] plot(\x,{(--1.6-0.0*\x)/1.0});
\draw [line width=1.2000000000000002pt,color=cyan,domain=-4.3:18.7] plot(\x,{(--0.2-0.0*\x)/1.0});
\draw [line width=1.2000000000000002pt,color=ffqqqq,domain=-4.3:18.7] plot(\x,{(-0.5-0.0*\x)/1.0});
\draw [line width=1.2000000000000002pt,color=ffqqqq,domain=-4.3:18.7] plot(\x,{(-0.65-0.0*\x)/1.0});
\draw [line width=1.2000000000000002pt,color=cyan] (3.677124344467705,-0.5)-- (6.322875655532295,-0.5);
\draw [line width=1.2000000000000002pt,color=cyan] (9.62,-0.5)-- (10.52,-0.5);
\draw [line width=1.2000000000000002pt,color=qqffqq] (-2.24,0.2)-- (9.48,0.2);
\end{scope}
\begin{scope}
\clip(10.0,-0.0) circle (0.848528137423857cm);
\draw (-1.2000000000000002,-5.320000000000002) -- (-1.2000000000000002,6.3);
\draw (-0.8000000000000002,-5.320000000000002) -- (-0.8000000000000002,6.3);
\draw (-0.40000000000000013,-5.320000000000002) -- (-0.40000000000000013,6.3);
\draw (-0.0,-5.320000000000002) -- (-0.0,6.3);
\draw (0.4,-5.320000000000002) -- (0.4,6.3);
\draw (0.8,-5.320000000000002) -- (0.8,6.3);
\draw (1.2000000000000002,-5.320000000000002) -- (1.2000000000000002,6.3);
\draw (1.6,-5.320000000000002) -- (1.6,6.3);
\draw (2.0,-5.320000000000002) -- (2.0,6.3);
\draw (2.4,-5.320000000000002) -- (2.4,6.3);
\draw (2.8,-5.320000000000002) -- (2.8,6.3);
\draw (3.1999999999999997,-5.320000000000002) -- (3.1999999999999997,6.3);
\draw (3.5999999999999996,-5.320000000000002) -- (3.5999999999999996,6.3);
\draw (4.0,-5.320000000000002) -- (4.0,6.3);
\draw (4.4,-5.320000000000002) -- (4.4,6.3);
\draw (4.800000000000001,-5.320000000000002) -- (4.800000000000001,6.3);
\draw (5.200000000000001,-5.320000000000002) -- (5.200000000000001,6.3);
\draw (5.600000000000001,-5.320000000000002) -- (5.600000000000001,6.3);
\draw (6.0,-5.320000000000002) -- (6.0,6.3);
\draw (6.4,-5.320000000000002) -- (6.4,6.3);
\draw (6.800000000000001,-5.320000000000002) -- (6.800000000000001,6.3);
\draw (7.200000000000001,-5.320000000000002) -- (7.200000000000001,6.3);
\draw (7.600000000000001,-5.320000000000002) -- (7.600000000000001,6.3);
\draw (8.0,-5.320000000000002) -- (8.0,6.3);
\draw (8.4,-5.320000000000002) -- (8.4,6.3);
\draw (8.8,-5.320000000000002) -- (8.8,6.3);
\draw (9.200000000000001,-5.320000000000002) -- (9.200000000000001,6.3);
\draw (9.600000000000001,-5.320000000000002) -- (9.600000000000001,6.3);
\draw (10.0,-5.320000000000002) -- (10.0,6.3);
\draw (10.4,-5.320000000000002) -- (10.4,6.3);
\draw (10.8,-5.320000000000002) -- (10.8,6.3);
\draw (11.200000000000001,-5.320000000000002) -- (11.200000000000001,6.3);
\draw (11.600000000000001,-5.320000000000002) -- (11.600000000000001,6.3);
\draw (12.0,-5.320000000000002) -- (12.0,6.3);
\draw [domain=-4.3:18.7] plot(\x,{(-2.0-0.0*\x)/1.0});
\draw [domain=-4.3:18.7] plot(\x,{(-1.4-0.0*\x)/1.0});
\draw [domain=-4.3:18.7] plot(\x,{(-0.7999999999999999-0.0*\x)/1.0});
\draw [domain=-4.3:18.7] plot(\x,{(-0.19999999999999996-0.0*\x)/1.0});
\draw [domain=-4.3:18.7] plot(\x,{(--0.4-0.0*\x)/1.0});
\draw [domain=-4.3:18.7] plot(\x,{(--1.0-0.0*\x)/1.0});
\draw [domain=-4.3:18.7] plot(\x,{(--1.6-0.0*\x)/1.0});
\draw [line width=1.2000000000000002pt,color=cyan,domain=-4.3:18.7] plot(\x,{(--0.2-0.0*\x)/1.0});
\draw [line width=1.2000000000000002pt,color=ffqqqq,domain=-4.3:18.7] plot(\x,{(-0.5-0.0*\x)/1.0});
\draw [line width=1.2000000000000002pt,color=ffqqqq,domain=-4.3:18.7] plot(\x,{(-0.65-0.0*\x)/1.0});
\draw [line width=1.2000000000000002pt,color=cyan] (3.677124344467705,-0.5)-- (6.322875655532295,-0.5);
\draw [line width=1.2000000000000002pt,color=cyan] (9.62,-0.5)-- (10.52,-0.5);
\draw [line width=1.2000000000000002pt,color=qqffqq] (-2.24,0.2)-- (9.48,0.2);
\end{scope}
\end{tikzpicture}
\caption{Topological argument to choose $T$ in the proof of lemma \ref{exclusive}}
\label{trivcov}
\end{figure}
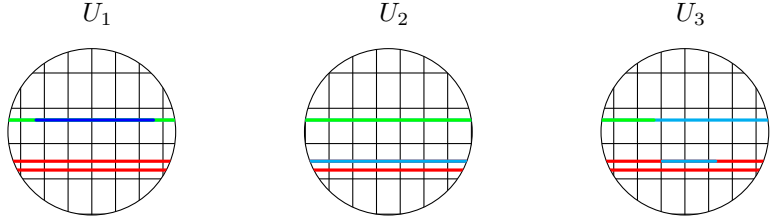

This argument is depicted on figure \ref{trivcov}, with $m_p([-1,1])$ in \textcolor{blue}{blue}, the remaining part of $m_p([-2,2])$ in \textcolor{green}{green}, the remaining part of $m_p([-3,3])$ in \textcolor{cyan}{cyan}, and the remaining part of $m_p(\R)$ in \textcolor{red}{red}. On this figure, $T=3$ is fine.

$U$ is open in $X$: let $x\in U$ and $t_0\in\R$ such that $l_x(t_0)=s_0\in S$. Then, $y\in X\mapsto l_y(t_0)\in X$ being a local diffeomorphism, we restrict it to a diffeomorphism between neighbourhoods $V_x$ and $V_{s_0}$ on which the north-east grid is trivial. Since $U$ contains a neighbourhood of $S$, $V_{s_0}\subset U$. Composing by $l_{\cdot}(-t_0)$ yields $V_x\subset U$.

And $U$ is also closed: let $x\in X$ approached by a sequence $(x_n)\in U^{\N}$. We have $(t_n),(s_n)$ such that $l_{x_n}(t_n)=s_n\in S$. Since $S$ is compact, we can assume that both $(s_n,t_n)$ converge to $(s,t)\in S\times\R\cup\{-\infty,+\infty\}$. If $x\in l_s(\R)\subset U$, it's fine. Else, $x$ can actually be approached by a sequence in $l_s(\R)$; in a neighbourhood of $x$ trivialising the north-east grid, $l_s(\R)$ is just a set of parallel lines of northern direction. So there are some $a,b\in\R$ such that $m_x(a)=l_s(b)$. Since north and east define locally constant vector fields with zero Lie bracket, their flow commute, so $s'=m_{s}(-a)\in S$ satisfies $l_{s'}(b)=x$.

By connectedness, $U=X$.}
 
For now on, let $S$ be as in the previous lemma.

\begin{figure}
\centering
\definecolor{ffwwqq}{rgb}{1.0,0.4,0.0}
\definecolor{qqccqq}{rgb}{0.0,0.8,0.0}
\definecolor{qqzzcc}{rgb}{0.0,0.6,0.8}
\definecolor{ffqqtt}{rgb}{1.0,0.0,0.2}
\begin{tikzpicture}[line cap=round,line join=round,>=triangle 45,scale=0.65]
\clip(-3.813333333333334,-1.846666666666667) rectangle (6.000000000000001,5.593333333333331);
\fill[line width=0.0pt,color=ffqqtt,fill=ffqqtt,fill opacity=0.5] (-3.0,1.0) -- (-0.41999999999999993,1.0) -- (-0.41999999999999993,6.3) -- (-3.0,6.3) -- cycle;
\fill[line width=0.0pt,color=ffqqtt,fill=ffqqtt,fill opacity=0.5] (4.0,1.0) -- (1.42,1.0) -- (1.42,-2.44) -- (4.0,-2.44) -- cycle;
\fill[line width=0.0pt,color=qqzzcc,fill=qqzzcc,fill opacity=0.5] (-0.41999999999999993,1.0) -- (1.7000000000000002,1.0) -- (1.7000000000000002,6.3) -- (-0.41999999999999993,6.3) -- cycle;
\fill[line width=0.0pt,color=qqccqq,fill=qqccqq,fill opacity=0.5] (2.6279999999999992,1.0) -- (2.6279999999999992,6.3) -- (1.7000000000000002,6.3) -- (1.7000000000000002,1.0) -- cycle;
\fill[line width=0.0pt,color=ffwwqq,fill=ffwwqq,fill opacity=0.5] (4.0,1.0) -- (4.0,6.3) -- (2.6279999999999992,6.3) -- (2.6279999999999992,1.0) -- cycle;
\fill[line width=0.0pt,color=qqzzcc,fill=qqzzcc,fill opacity=0.5] (-0.8799999999999999,1.0) -- (-3.0,1.0) -- (-3.0,-2.44) -- (-0.8799999999999999,-2.44) -- cycle;
\fill[line width=0.0pt,color=ffwwqq,fill=ffwwqq,fill opacity=0.5] (0.4920000000000009,1.0) -- (-0.8799999999999999,1.0) -- (-0.8799999999999999,-2.44) -- (0.4920000000000009,-2.44) -- cycle;
\fill[line width=0.0pt,color=qqccqq,fill=qqccqq,fill opacity=0.5] (1.42,1.0) -- (0.4920000000000009,1.0) -- (0.4920000000000009,-2.44) -- (1.42,-2.44) -- cycle;
\draw [line width=1.6pt] (-3.0,1.0)-- (4.0,1.0);
\draw (-3.0,-1.846666666666667) -- (-3.0,5.593333333333331);
\draw (4.0,-1.846666666666667) -- (4.0,5.593333333333331);
\draw (-0.41999999999999993,1.0) -- (-0.41999999999999993,5.593333333333331);
\draw (1.7000000000000002,1.0) -- (1.7000000000000002,5.593333333333331);
\draw (2.6279999999999992,1.0) -- (2.6279999999999992,5.593333333333331);
\draw (1.42,1.0) -- (1.42,-1.846666666666667);
\draw (-0.8799999999999999,1.0) -- (-0.8799999999999999,-1.846666666666667);
\draw (0.4920000000000009,1.0) -- (0.4920000000000009,-1.846666666666667);
\begin{scriptsize}
\draw [fill=black] (-0.4199999999999999,3.0733333333333315) circle (1.5pt);
\draw [fill=black] (1.7000000000000004,5.1) circle (1.5pt);
\draw [fill=black] (-0.8799999999999999,-0.08666666666666739) circle (1.5pt);
\end{scriptsize}
\end{tikzpicture}
\caption{Strip decomposition locally around $S_0$ (black points are singularities)}
\label{localdraw}
\end{figure}

\demode{of Proposition \ref{flatRS}}{Let $\eps>0$ such that $l\,:\,(s,t)\in S\times(-\eps,\eps)\mapsto l_s(t)\in X$ is injective. Let $U$ be its image. Since $U\subset X$ has non zero Lebesgue measure, by Poincaré recurrence theorem, for almost all $s\in S$ there is a time $t_U(s)>\eps$ such that $l_s(t_U(s))\in U$ and hence $t(s)\in (t_U(s)-\eps,t_U(s)+\eps)$ such that $l_s(t(s))\in S$. Without loss of generality, we can assume that $t(s)>0$ is minimal and call it the \emph{first return time} of $l_s$. By \emph{almost all}, we understand $S_1\subset S$ such that $l(S_1\times (-\eps,\eps))$ and $l(S\times (-\eps,\eps))$ have the same Lebesgue measure.

Let us cover $X$ by a finite set ${\cal C}$ of open balls trivialising the north-east grid. If $x\in S_1$ and $l_x$ does not go through any singular point before its first return to $S$, then there is a neighbourhood $x\in V_x\subset S$ such that for all $t\in [0,t(x)]$ and for all $y\in V_x$, an open ball of ${\cal C}$ contains both $l_y(t)$ and $l_x(t)$. Then all $l_y$ for $y\in V_x$ are parallel paths with constant distance to each other before their first return to $S$ and they all have the same first return time $t(y)=t(x)$. Particularly, $V_x\subset S_1$. So $S=S_1$ by connectedness. 

Hence, this construction defines a striped decomposition of $X$ in the following sense (see figure \ref{localdraw}): a \emph{strip} is an element of $$\{l_x(t)\mid x\in S,t\in [0,t(x)),l_x(s)\mbox{ is not a singular point for }s< t(x)\}\;/\hspace{-2pt}\sim$$
\noindent where $\sim$ is the finest equivalence relation such that $l_x(s)\sim l_x(t)$ if $x\in S,s,t\in [0,t(x))$ and $x\sim y$ if $x\in S$ and $y\in V_x$. According to Lemma \ref{exclusive}, strips partition $X$. This gives equivalently a construction of $X$ from glued rectangles as in figure \ref{rectangles}.

\begin{figure}
\centering
\definecolor{ffwwqq}{rgb}{1.0,0.4,0.0}
\definecolor{qqccqq}{rgb}{0.0,0.8,0.0}
\definecolor{qqzzcc}{rgb}{0.0,0.6,0.8}
\definecolor{ffqqtt}{rgb}{1.0,0.0,0.2}
\begin{tikzpicture}[line cap=round,line join=round,>=triangle 45,x=0.9cm,y=0.9cm]
\clip(-3.486000000000001,0.5146666666666685) rectangle (4.668666666666672,6.718666666666666);
\draw [line width=1.6pt,color=ffqqtt] (-3.0,6.3)-- (-0.41999999999999993,6.3);
\draw [line width=1.6pt,color=qqzzcc] (-0.41999999999999993,3.6)-- (1.7000000000000002,3.6);
\draw [line width=1.6pt,color=qqccqq] (1.7000000000000002,4.74)-- (2.6279999999999992,4.74);
\draw [line width=1.6pt,color=ffwwqq] (2.6279999999999992,5.94)-- (4.0,5.94);
\draw [line width=1.6pt,color=qqzzcc] (-3.0,1.0)-- (-0.8799999999999999,1.0);
\draw [line width=1.6pt,color=ffwwqq] (-0.8799999999999999,1.0)-- (0.4920000000000009,1.0);
\draw [line width=1.6pt,color=qqccqq] (0.4920000000000009,1.0)-- (1.4200000000000002,1.0);
\draw [line width=1.6pt,color=ffqqtt] (1.4200000000000002,1.0)-- (4.0,1.0);
\draw (1.7000000000000002,4.74)-- (1.7000000000000002,3.6);
\draw (1.7660000000000002,4.195666666666667) -- (1.6340000000000001,4.195666666666667);
\draw (1.7660000000000002,4.144333333333334) -- (1.6340000000000001,4.144333333333334);
\draw (-3.0,1.0)-- (-3.0,2.499999999999999);
\draw (-3.0733333333333337,1.75) -- (-2.9266666666666667,1.75);
\draw (-3.0,2.499999999999999)-- (-3.0,6.3);
\draw (-3.0,4.553999999999999) -- (-2.9010000000000002,4.476999999999999);
\draw (-3.0,4.553999999999999) -- (-3.099,4.476999999999999);
\draw (-3.0,4.3999999999999995) -- (-2.9010000000000002,4.3229999999999995);
\draw (-3.0,4.3999999999999995) -- (-3.099,4.3229999999999995);
\draw (-0.41999999999999993,3.6)-- (-0.41999999999999993,4.800000000000001);
\draw (-0.4200000000000002,4.277) -- (-0.3210000000000001,4.2);
\draw (-0.4200000000000002,4.277) -- (-0.5190000000000002,4.2);
\draw (-0.41999999999999993,4.800000000000001)-- (-0.41999999999999993,6.3);
\draw (-0.49333333333333357,5.55) -- (-0.3466666666666668,5.55);
\draw (4.0,1.0)-- (4.0,2.14);
\draw (3.9339999999999993,1.5443333333333336) -- (4.066,1.5443333333333336);
\draw (3.9339999999999993,1.595666666666667) -- (4.066,1.595666666666667);
\draw (4.0,2.14)-- (4.0,5.94);
\draw (3.9999999999999996,4.194) -- (4.098999999999999,4.117);
\draw (3.9999999999999996,4.194) -- (3.9009999999999994,4.117);
\draw (3.9999999999999996,4.04) -- (4.098999999999999,3.9630000000000005);
\draw (3.9999999999999996,4.04) -- (3.9009999999999994,3.9630000000000005);
\draw (2.6279999999999992,4.74)-- (2.6279999999999992,5.94);
\draw (2.627999999999999,5.417000000000001) -- (2.726999999999999,5.340000000000001);
\draw (2.627999999999999,5.417000000000001) -- (2.5289999999999986,5.340000000000001);
\begin{scriptsize}
\draw [fill=black] (-0.41999999999999993,6.3) circle (1.5pt);
\draw [fill=black] (-3.0,6.3) circle (1.5pt);
\draw [fill=black] (-0.41999999999999993,3.6) circle (1.5pt);
\draw [fill=black] (1.7000000000000002,3.6) circle (1.5pt);
\draw [fill=black] (1.7000000000000002,4.74) circle (1.5pt);
\draw [fill=black] (2.6279999999999992,4.74) circle (1.5pt);
\draw [fill=black] (2.6279999999999992,5.94) circle (1.5pt);
\draw [fill=black] (4.0,5.94) circle (1.5pt);
\draw [fill=black] (-3.0,1.0) circle (1.5pt);
\draw [fill=black] (-0.8799999999999999,1.0) circle (1.5pt);
\draw [fill=black] (0.4920000000000009,1.0) circle (1.5pt);
\draw [fill=black] (1.4200000000000002,1.0) circle (1.5pt);
\draw [fill=black] (4.0,1.0) circle (1.5pt);
\draw [fill=black] (-0.41999999999999993,4.800000000000001) circle (1.5pt);
\draw [fill=black] (-3.0,2.499999999999999) circle (1.5pt);
\draw [fill=black] (4.0,2.14) circle (1.5pt);
\end{scriptsize}
\end{tikzpicture}
\caption{Gluing rectangles to construct $X$}
\label{rectangles}
\end{figure}

Hence, $X$ is a $PK^0$ curve.  
}

\subsection{Does a $PK^0$ curve cover of a torus?}

Let us come to the proof of theorem A. First, we state what we mean by $V$ and $H$-periods.

\defi{Let $X$ be a $PK^0$ surface with a complex flat structure with associated 1-form $\de z$; $X$ is obtained from gluing rectangles as in the proof of Proposition \ref{flatRS}. For any closed path $\gamma$ on $X$,
$$\intg{\gamma}{}\de z = a+\i b\in \C,$$
\noindent and $a$ (respectively $b$) is called an \emph{$H$-period} (respectively a \emph{$V$-period}) of $X$. We may refer to $a$ as the $H$-period on the path $\gamma$.}

\rem[finitegener]{The set of all $H$-periods of $X$ is a subgroup $P_H$ of $\R$. Actually, it is finitely generated, by the $H$-periods on horizontal paths between singular points of $X$ (e.g. on the figure \ref{rectangles}, these paths are only the coloured segments).}

\prop[notcover]{If both $P_H$ and $P_V$ are discrete subgroups of $\R$ (or equivalently the set of all periods is a discrete subgroup of $\C$), then $X$ is a ramified covering of the complex torus $\C/\Gamma$, where $\Gamma=P_H+\i P_V$.}

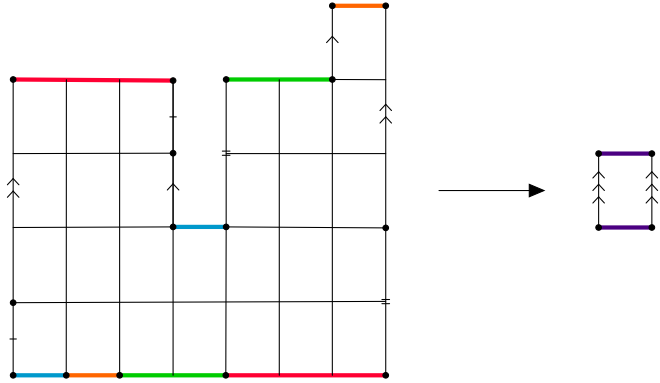
\begin{figure}
\centering
\definecolor{ubqqys}{rgb}{0.29411764705882354,0.0,0.5098039215686274}
\definecolor{ffwwqq}{rgb}{1.0,0.4,0.0}
\definecolor{qqccqq}{rgb}{0.0,0.8,0.0}
\definecolor{qqzzcc}{rgb}{0.0,0.6,0.8}
\definecolor{ffqqtt}{rgb}{1.0,0.0,0.2}
\begin{tikzpicture}[line cap=round,line join=round,>=triangle 45,x=1cm,y=1.4cm,scale=0.7]
\clip(-3.8262666666666667,0.24773333333333483) rectangle (9.709600000000012,6.765600000000007);
\draw [line width=1.6pt,color=ffqqtt] (-3.0,5.002666666666666)-- (0.004666666666669261,4.9879999999999995);
\draw [line width=1.6pt,color=qqzzcc] (0.004666666666669261,3.008)-- (1.0019999999999998,3.008);
\draw [line width=1.6pt,color=qqccqq] (1.0019999999999998,5.002666666666666)-- (2.9966666666666733,5.002666666666666);
\draw [line width=1.6pt,color=ffwwqq] (2.9966666666666733,6.0)-- (4.0,6.0);
\draw [line width=1.6pt,color=qqzzcc] (-3.0,1.0)-- (-2.0026666666666695,1.0);
\draw [line width=1.6pt,color=ffwwqq] (-2.0026666666666695,1.0)-- (-0.9993333333333427,1.0);
\draw [line width=1.6pt,color=qqccqq] (-0.9993333333333427,1.0)-- (0.9953333333333306,1.0);
\draw [line width=1.6pt,color=ffqqtt] (0.9953333333333306,1.0)-- (4.0,1.0);
\draw (1.0019999999999998,5.002666666666666)-- (1.0019999999999998,3.008);
\draw (1.0746,4.033566666666666) -- (0.9293999999999999,4.033566666666666);
\draw (1.0746,3.9771) -- (0.9293999999999999,3.9771);
\draw (-3.0,1.0)-- (-3.0,1.9826666666666657);
\draw (-3.0605,1.491333333333333) -- (-2.9395,1.491333333333333);
\draw (-3.0,1.9826666666666657)-- (-3.0,5.002666666666666);
\draw (-3.0,3.6620666666666657) -- (-2.8911,3.577366666666666);
\draw (-3.0,3.6620666666666657) -- (-3.1089,3.577366666666666);
\draw (-3.0,3.4926666666666657) -- (-2.8911,3.4079666666666655);
\draw (-3.0,3.4926666666666657) -- (-3.1089,3.4079666666666655);
\draw (0.004666666666669261,3.008)-- (0.004666666666669261,4.005333333333334);
\draw (0.004666666666669596,3.5913666666666675) -- (0.11356666666666969,3.5066666666666673);
\draw (0.004666666666669596,3.5913666666666675) -- (-0.1042333333333305,3.5066666666666673);
\draw (0.004666666666669261,4.005333333333334)-- (0.004666666666669261,4.9879999999999995);
\draw (-0.05583333333333045,4.496666666666666) -- (0.06516666666666965,4.496666666666666);
\draw (4.0,1.0)-- (4.0,2.9946666666666664);
\draw (3.9274000000000004,1.9690999999999992) -- (4.0726,1.9690999999999992);
\draw (3.9274000000000004,2.025566666666666) -- (4.0726,2.025566666666666);
\draw (4.0,2.9946666666666664)-- (4.0,6.0);
\draw (4.0,4.666733333333333) -- (4.1089,4.582033333333333);
\draw (4.0,4.666733333333333) -- (3.8911000000000002,4.582033333333333);
\draw (4.0,4.497333333333333) -- (4.1089,4.412633333333333);
\draw (4.0,4.497333333333333) -- (3.8911000000000002,4.412633333333333);
\draw (2.9966666666666733,5.002666666666666)-- (2.9966666666666733,6.0);
\draw (2.9966666666666737,5.586033333333334) -- (3.105566666666674,5.501333333333333);
\draw (2.9966666666666737,5.586033333333334) -- (2.8877666666666735,5.501333333333333);
\draw (-2.0026666666666695,1.0)-- (-2.0,5.0);
\draw (-0.9993333333333427,1.0)-- (-1.0,5.0);
\draw (0.004666666666669261,4.9879999999999995)-- (0.004666666666669261,1.0);
\draw (1.0019999999999998,3.008)-- (0.9953333333333306,1.0);
\draw (2.0,5.0)-- (2.0,1.0);
\draw (2.9966666666666733,5.002666666666666)-- (2.9966666666666733,1.0);
\draw (-3.0000000000000004,4.0)-- (0.004666666666669261,4.005333333333334);
\draw (-3.0000000000000004,3.0000000000000004)-- (0.004666666666669261,3.008);
\draw (-3.0,1.9826666666666657)-- (4.0,2.0);
\draw (1.0019999999999998,3.008)-- (4.0,2.9946666666666664);
\draw (1.0019999999999998,4.0)-- (4.0,4.0);
\draw (2.9966666666666733,5.002666666666666)-- (4.0,5.0);
\draw [line width=1.6pt,color=ubqqys] (8.0,4.0)-- (9.0,4.0);
\draw [line width=1.6pt,color=ubqqys] (9.0,3.0)-- (8.0,3.0);
\draw (8.0,3.0)-- (8.0,4.0);
\draw (8.0,3.5846999999999998) -- (8.1089,3.4999999999999996);
\draw (8.0,3.5846999999999998) -- (7.8911,3.4999999999999996);
\draw (8.0,3.4153) -- (8.1089,3.3305999999999996);
\draw (8.0,3.4153) -- (7.8911,3.3305999999999996);
\draw (8.0,3.7540999999999998) -- (8.1089,3.6694);
\draw (8.0,3.7540999999999998) -- (7.8911,3.6694);
\draw (9.0,3.0)-- (9.0,4.0);
\draw (9.000000000000002,3.5846999999999998) -- (9.108900000000002,3.4999999999999996);
\draw (9.000000000000002,3.5846999999999998) -- (8.891100000000002,3.4999999999999996);
\draw (9.000000000000002,3.4153) -- (9.108900000000002,3.3305999999999996);
\draw (9.000000000000002,3.4153) -- (8.891100000000002,3.3305999999999996);
\draw (9.000000000000002,3.7540999999999998) -- (9.108900000000002,3.6694);
\draw (9.000000000000002,3.7540999999999998) -- (8.891100000000002,3.6694);
\draw [->] (5.0,3.5) -- (7.0,3.5);
\begin{scriptsize}
\draw [fill=black] (0.004666666666669261,4.9879999999999995) circle (1.5pt);
\draw [fill=black] (-3.0,5.002666666666666) circle (1.5pt);
\draw [fill=black] (0.004666666666669261,3.008) circle (1.5pt);
\draw [fill=black] (1.0019999999999998,3.008) circle (1.5pt);
\draw [fill=black] (1.0019999999999998,5.002666666666666) circle (1.5pt);
\draw [fill=black] (2.9966666666666733,5.002666666666666) circle (1.5pt);
\draw [fill=black] (2.9966666666666733,6.0) circle (1.5pt);
\draw [fill=black] (4.0,6.0) circle (1.5pt);
\draw [fill=black] (-3.0,1.0) circle (1.5pt);
\draw [fill=black] (-2.0026666666666695,1.0) circle (1.5pt);
\draw [fill=black] (-0.9993333333333427,1.0) circle (1.5pt);
\draw [fill=black] (0.9953333333333306,1.0) circle (1.5pt);
\draw [fill=black] (4.0,1.0) circle (1.5pt);
\draw [fill=black] (0.004666666666669261,4.005333333333334) circle (1.5pt);
\draw [fill=black] (-3.0,1.9826666666666657) circle (1.5pt);
\draw [fill=black] (4.0,2.9946666666666664) circle (1.5pt);
\draw [fill=black] (9.0,3.0) circle (1.5pt);
\draw [fill=black] (8.0,3.0) circle (1.5pt);
\draw [fill=black] (9.0,4.0) circle (1.5pt);
\draw [fill=black] (8.0,4.0) circle (1.5pt);
\end{scriptsize}
\end{tikzpicture}
\caption{The finite sheeted ramified covering $\tilde{X}\to\C/\Gamma$}
\label{covering}
\end{figure}

\demo{Let $\Gamma=P_H+\i P_V$. Consider the rectangle decomposition of $X$ as a (simply-connected) polygonal closed set $\tilde{X}\subset \C$ from which $X$ is glued by the projection $p\,:\,\tilde{X}\to X$. Fixing a base point and integrating $\de z$ along paths contained in $\tilde{X}$, we define an injective, holomorphic map: $\tilde{X}\to\C$. Composed with the canonical $\C\to\C/\Gamma$, it yields a finite-sheeted surjective map $\tilde{X}\to\C/\Gamma$ which is a covering inside $\tilde{X}$. It is depicted in figure \ref{covering}. 

By definition of $\Gamma$, this factorizes through $p$ and we get a map: $X\to\C/\Gamma$, which is a finite-sheeted covering ramified only at points corresponding to vertices of $\tilde{X}$, that is to say at the singular points of $X$.}

We are now left with the converse statement of theorem A.

\prop[periodsarediscrete]{Let $X$ be a $PK^0$ curve and a ramified covering of a torus $\pi\,:\, X\to\C/\Gamma$ as well. Then $P_H$ and $P_V$ both are discrete subgroups of $\R$.}

\demo{It is equivalent to show that $P_H+\i P_V$ is a discrete subgroup of $\C$, which is clear since:
\begin{align*}
P_H+\i P_V
=\left\lbrace\intg{\gamma}{}\de z\mid\gamma\mbox{ closed path on }X\right\rbrace
\subset\left\lbrace\intg{\gamma'}{}\de z\mid\gamma'\mbox{ closed path on }\C/\Gamma\right\rbrace
=\Gamma.
\end{align*}}

\rem{By Remark \ref{finitegener}, if $X$ is given with its rectangle decomposition, we can easily construct a finite number of generators of $P_H$ and $P_V$. Checking that $P_H$ and $P_V$ are discrete is then tantamount to checking that a few generators are commensurable, which is quite an effective criterion.}

\subsection{Generalising to products of $PK^0$ curves}

Let us now prove theorem B.

\demo{The converse direction is clear.

Let now $\pi\,:\, X=C_1\times C_2\to\T$ be a ramified covering of a torus by a product of $PK^0$ curves. Let us prove that both curves $C_1$ and $C_2$ are ramified coverings of complex one-dimensional tori.

Denote by $\de z_1,\de z_2$ the unique (up to a complex multiplicative constant) holomorphic flat differentials on $C_1,C_2$, and take a basis of flat differentials $\de z,\de w$ on $\T$. For any fixed $b\in C_2$, we pull them back by
$$f_b\,:\, C_1\simeq C_1\times\{b\}\hookrightarrow C_1\times C_2\overset{\pi}{\longrightarrow}\T.$$

These pull-backs are holomorphic flat differentials on $C_1$, so multiples of $\de z_1$, and their periods are the same as the periods of the forms $\de z,\de w$ on the torus, so discrete subgroups of $\C$. But by Proposition \ref{notcover}, the periods of a non-zero multiple of $\de z_1$ form a discrete subgroup of $C_1$ if and only if $C_1$ is a ramified covering of a torus. 

Hence, if we assume that $C_1$ is not a ramified covering of a torus (working by contradiction), we get $f_b^*(\de z)=f_b^*(\de w)=0$, that is to say that at a point $(a,b)\in C_1\times C_2$, both linear forms $(\pi^*\de z)_{(a,b)}$ and $(\pi^*\de w)_{(a,b)}$ on the tangent space $T_{(a,b)}(X)$ have the same one-dimensional kernel $T_{(a,b)}(C_1\times\{b\})$. So these two linear forms $\pi^*\de z,\pi^*\de w$ are proportional. But $\de z,\de w$ were set linearly independent on $\T$, absurd! 

So $C_1$ is a ramified covering of a torus.
}

\rem{The compatibility of our pull-backs of holomorphic 1-forms with the flat structures arises from the fact that we always pull back through local isometries (for $\pi$, it is example 3.4.4 in \cite{BBIv}). It is very important since our main criterion, that is theorem A, focuses on periods for a holomorphic 1-form $\de z$ compatible with the flat structure.}

\section{Characterizing $PK^0$ product of curves}

\subsection{Some more definitions}

Let $X$ be a $PK^0$ surface with 2-skeleton $X_s$. 

\defi{A \emph{complex direction} of $X$ is defined as a parallel line bundle on $X\setminus X_s$ which is also a subsheaf of $\Omega^1(X\setminus X_s)$.}

Let $F$ be a codimension 2 face of $X$ and $F_s$ be the set of points $x\in F$ which are codimension 4 singularities of $X$. In a neighbourhood of a point $x\in F\setminus F_s$, there is a chart $(z,z_1)\in U_x\subset X\mapsto (z^n,z_1)\in\C^2$ such that locally $F\cap U_x=\{z=0\}$. 

\defi{\emph{The complex direction of $F$} is the unique complex direction $L_F$ of $X$ defined on $X\setminus X_s\cup F\setminus F_s$ such that $L_F(U_x)=\C\cdot\de z\subset\Omega^1(X\setminus X_s\cup F\setminus F_s)$. }

\defi{A complex direction $L$ of $X$ is said \emph{relevant} if there is a codimension 2 face $F$ of $X$ such that $L=L_F$.}

Examining whether $X$ could be a product of curves, we will be interested in eventual splittings of its cotangent bundle. Our first goal is to extend certain complex directions of $X$ to the whole surface. This works out when two or less distinct complex directions of $X$ are relevant.

Then we give an example of a $PK^0$ surface $X$ with three relevant complex directions, which is not a product of two curves. It arises as a ramified covering of a torus over some elliptic curves. Let us summarize some useful facts about flat elliptic curves inside a complex 2-dimensional torus, and deduce some easy statements on ramified coverings of tori and $PK^0$ surfaces.

\subsection{Excursion in the realm of complex tori}

\defi{Two complex tori $\T_1,\T_2$ of same complex dimension $n$ are said \emph{isogenous} if one of these equivalent conditions is satisfied:
\begin{itemize}
\item there is a finite-sheeted covering $\T_1\to \T_2$ ;
\item there is a finite-sheeted covering $\T_2\to \T_1$ ;
\item there are lattices $\Lambda_1,\Lambda_2$ of $\C^n$ such that $\T_i$ is biholomorphic to $\C^n/\Lambda_i$ and for some non-zero integers $n_1,n_2$, $n_1\Lambda_1=n_2\Lambda_2$.
\end{itemize}
Isogeny defines an equivalence relation on the moduli space of complex $n$-dimensional tori: we denote it ambiguously by $\simeq $.}

\defi{Given a flat elliptic curve $C$ inside a complex $n$-dimensional torus $\T$, there is a holomorphic local flat coordinate $z$ on $\T$ such that $C$ is locally given by an equation of the form $z=c$ for a constant $c\in\C^n$. This coordinate $z$ is actually globally defined up to an additive constant. Let us call it \emph{the coordinate induced by $C$}: up to additive and multiplicative constant, it is unique. 

Two flat elliptic curves $C_1,C_2$ in $\T$ are said \emph{parallel} if their induced coordinates can be made equal by constant shift and rescaling.}

\lem[prod]{Let $\T$ be a flat two-dimensional complex torus containing two non parallel flat elliptic curves $C_1,C_2$. Then $\T\simeq C_1\times C_2$.}

\demo{Let $z_1,z_2$ be the coordinates of $\T$ induced by $C_1,C_2$. Since they are globally defined up to additive constant, they induce global holomorphic maps $\tilde{z_1}\,:\, \T\to\C/\Gamma_1,\tilde{z_2}\,:\,\T\to\C/\Gamma_2$ for some lattices $\Gamma_1,\Gamma_2$. The equation $\tilde{z_i}=c_i$ is now a global equation for $C_i$ (for a good choice of $c_i$) and there is a projection in the local basis $\tilde{z_1},\tilde{z_2}$ given by $(x_1,x_2)\in\T\mapsto ((c_1,x_1),(c_2,x_2))\in C_1\times C_2$. Hence, $\T$ is a finite-sheeted covering of $C_1\times C_2$.}

We can now prove Proposition \ref{covtispk}.

\demode{of Proposition \ref{covtispk}}{Let $\T$ be a torus and $C_1,\ldots,C_n$ flat curves in $\T$ with no three of them having a common intersection point. Let $X$ be a ramified covering of $\T$ over $C_1,\ldots,C_n$. It has trivial holonomy; it is only left to check the local condition that it should be polyhedral Kähler.

It is clear outside of the ramification locus. In a neighbourhood of a point belonging to only one ramification line, there are two local coordinate $z_1,z_2$ as in the proof of Lemma \ref{prod} making the situation isometric to: $(z_1,z_2)\mapsto(z_1^m,z_2)\in\C^2$ over a neighbourhood of the point $(0,0)$, with ramification line $\{z_1=0\}$. Take any $\eps>0$. Gluing cubes $$\mathrm{Cub}_k=\{(z_1,z_2)\in\C^2\mid|z_1|,|z_2|<\eps,\, \mathrm{arg}(z_1)\in [\pi k,\pi (k +1)]\}$$ by isometries for $k\in\Z/ 2m\Z$, we get a local $PK^0$ model of our singularity. 

In a neighbourhood of an intersection point of two ramification lines, the situation is again as in Lemma \ref{prod}, that is to say like in the chart $(z_1,z_2)\mapsto(z_1^m,z_2^n)\in\C^2$ over a neighbourhood of $(0,0)$. Then consider all cubes $$\mathrm{Cub}_{k,l}=\{(z_1,z_2)\in\C^2\mid|z_1|,|z_2|<\eps,\, \mathrm{arg}(z_1)\in [\pi k,\pi (k +1)]\},\mathrm{arg}(z_2)\in [\pi l,\pi (l +1)]\}$$ for $k\in\Z/ 2m\Z,l\in\Z/ 2n\Z$ and glue $\mathrm{Cub}_{k,l}$ with $\mathrm{Cub}_{k+1,l}$ (respectively $\mathrm{Cub}_{k,l}$ with $\mathrm{Cub}_{k,l+1}$) on the side $\mathrm{arg}(z_1)=\pi (k+1)$ (respectively $\mathrm{arg}(z_2)=\pi (l+1)$) by the identity. This gives a local $PK^0$ model of our product-type singularity.}

\rem{Let $X$ be a $PK^0$ surface and a ramified covering of a complex flat two-dimensional torus $\pi\,:\, X\to\T$ over curves $C_1,\ldots,C_n\subset\T$. Let $X_s$ be as usual the singular locus of $X$ as a $PK^0$ surface. Then $\pi(X_s)$ is a union of flat elliptic curves inside $\T$, since any codimension 2 face of $X$ has a fixed complex direction. Conversely, each $\pi^{-1}(C_i)$ is a curve of fixed complex direction in $X$. There is generally no equality between $\pi(X_s)$ and $C_1\cup\ldots\cup C_n$.}

\rem{If a flat torus contains three flat elliptic curves $C_1,C_2,C_3$, none of them parallel, then the projections $C_3\to C_1,C_3\to C_2$ are non-zero, so $C_1,C_2,C_3$ are isogenous. To that extent, considering $PK^0$ surfaces being ramified over at least three curves on a torus is a very nice and specific situation.}

\subsection{Extend line bundles}

In the forthcoming proofs, we need to extend line bundles near some singularities. Let us prepare for it with the following results.

\defi[spt]{If $L$ is a line bundle on $X\setminus X_s$ and a subsheaf of $\Omega^1(X\setminus X_s)$, we will call a line bundle $\tilde{L}$ an \emph{extension of $L$ to the set $U$} if $X\setminus X_s\subset U\subset X$, $\tilde{L}$ is a subsheaf of $\Omega^1(U)$, and the pull-back of $\tilde{L}$ by the injection $X\setminus X_s\hookrightarrow U$ is equal to $L$.}

In order to extend line bundles near singularities of complex codimension 2, there is the lemma 3.2 from \cite{Cata} by Fabricio Catanese and Mateo Franciosi: 

\defi[ext]{Let $X$ be a complex surface. A \emph{special tensor} $\omega$ is a never-vanishing global section of the sheaf $$\mathrm{SpT}:=\mathrm{Sym}^2(\Omega^1(X))\otimes(\det\,\Omega^1(X))^{-1}\simeq\mathrm{End}^0(TX),$$
where $\mathrm{End}^0(TX)$ is the sheaf of trace zero endomorphisms of $TX$.
Locally, in a basis $\left(\frac{\partial}{\partial z_1},\frac{\partial}{\partial z_2}\right)$ of the tangent space, the isomorphism of sheaves is given by $$a_{11}\frac{\de z_1\otimes\de z_1}{\de z_1\wedge\de z_2}+a\frac{\de z_1\otimes\de z_2+\de z_2\otimes\de z_1}{\de z_1\wedge\de z_2}+a_{22}\frac{\de z_2\otimes\de z_2}{\de z_1\wedge\de z_2}\,\longmapsto\,
\left(\begin{matrix}
-a & a_{22} \\ 
a_{11} & a
\end{matrix}\right) .$$}

\rem{On a compact complex surface, the fact that a special tensor is never-vanishing is tantamount to the fact that its (holomorphic on $X$, hence constant) determinant as endomorphism of $TX$ is non-zero.}

\prop[isspt]{Let $X$ be a complex surface. A splitting of $\Omega^1(X)$ as direct sum of two line bundles $M_1,M_2$ is equivalent to the existence of a special tensor $\omega$ on $X$. As it happens, one can assume that locally, if $M_i=\C\cdot\de z_i$, then $\omega$ corresponds to the diagonal endomorphism of the tangent space with eigenvalues $-1,1$ and respective eigenvectors $\frac{\partial}{\partial z_1},\frac{\partial}{\partial z_2}$.}

Let us use theses results to prove:

\lem[hartogs]{Let $X$ is a complex surface with a finite set of marked points $S$. Suppose either that $X$ is compact or that each $s\in S$ is contained in a compact complex curve $C\subset X$. If there are two line bundles $M_1,M_2$ such that $\Omega^1(X\setminus S)=M_1\oplus M_2$, then one can extend them (in the sense of Definition \ref{ext}) as $\tilde{M_1},\tilde{M_2}$ to $X$. It yields a splitting $\Omega^1(X)=\tilde{M_1}\oplus\tilde{M_2}$.}

\demo{Combine Proposition \ref{spt} on the complex surface $X\setminus S$, the fact that the sheaf $\mathrm{SpT}$ is locally the same as a 3-dimensional free sheaf of holomorphic functions and Hartogs' extension theorem to define a global section of $\mathrm{SpT}$ on $X$. Its determinant is non-zero on $X\setminus S$; at a point $s\in S$ belonging to a compact complex curve $C$ (respectively to $X$, if it is compact), since the determinant should be holomorphic, hence constant and non-zero on $C$ (respectively $X$), we get $\det(\omega_s)\ne 0$, so $\omega_s\ne 0$. So we have a special tensor on $X$, and the converse direction of Proposition \ref{isspt} applies to split the cotangent bundle of $X$.}

Near complex codimension 1 singularities, we extend line bundles thanks to the following result:

\lem[localsplit]{Let $X$ be a $PK^0$ surface with 2-skeleton $X_s$ and $C\subset X_s$ an irreducible component. Let $C_s$ be the set of complex codimension 2 singularities of $X$ belonging to $C$. Let $x\in C\setminus C_s$ with a neighbourhood $V_x\subset X$ not intersecting $C_s$.
Let $L_1,L_2$ two holomorphic line bundles defined on $V_x\setminus C$ such that $L_1$ is the complex direction of $C$. Then $L_1,L_2$ can be extended to holomorphic line bundles on a neighbourhood $x\in U_x\subset V_x$.}

\demo{Reducing to a smaller neighbourhood $U_x\subset V_x$, we have a chart $\phi_x\,:\,(z_1,z_2)\in U_x\to (z_1^n,z_2)\in\C^2,$ in local coordinates $(z_1,z_2$). This coordinate system can be chosen such that $L_1,L_2$ are spanned by the 1-forms $\de z_1,\de z_2$ on $U_x\setminus C$.

Let us take local coordinates $z, z_2$ on $U_x$ pulled back from $(z_1,z_2)$ through this chart: $nz^{n-1}\de z=\de z_1,\de z_2=\de z_2$. Setting $L_1(U_x)=\C\cdot\de z$ and $L_2(U_x)=\C\cdot\de z_2$ makes sense and gives rise to the two extended holomorphic line bundles as wished.}

\rem{Complex codimension 1 singularities are tackled locally, whereas we need some global hypothesis (some compactness condition) to deal with complex codimension 2 singularities.}

\subsection{$PK^0$ surfaces with two complex direction}

We start by recalling a result of \cite{Beauv} to characterise a product of complex curves.

\theonm{(by Arnaud Beauville)}{Let $X$ be a compact complex surface. Its tangent bundle splits if and only if the universal cover of $X$ is a product of two complex surfaces $U\times V$ and the fundamental group $\pi_1(X)$ acts diagonally on $U\times V$, or $X$ is a Hopf surface with universal cover $\C^2\setminus\{0\}$ and $\pi_1(X)\simeq \Z\oplus\Z/m\Z$ generated by $(x,y)\mapsto (\alpha x,\beta y))$ with $|\alpha|\le|\beta|<1$ and $(x,y)\mapsto (\lambda x,\mu y)$ with $\lambda,\mu$ primitive $m$-roots of unity for some $m\in\N^*$.}

Since a Hopf surface can not be a $PK^0$ surface as proven in Section {\bf 1.2}, we are left to prove

\prop[noextra]{Let $X$ be a $PK^0$ surface with at most two complex directions. Then its cotangent bundle splits.}

\demo{Let $L_1,L_2$ be two locally independent line bundles defined on $X\setminus X_s$ by translating $\de z_1,\de z_2$ (thanks to flatness and trivial holonomy of $X$), such that the at most two complex directions of $X$ are among them.
 They extend to $\tilde{L_1},\tilde{L_2}$ on $X\setminus S$ by Lemma \ref{localsplit}, for $S\subset X_s$ the set of complex codimension 2 singularities on $X$; since $X$ is compact, Lemma \ref{hartogs} then extends $\tilde{L_1},\tilde{L_2}$ to $X$, and we are done.}

\subsection{Ramified torus over three curves of different complex direction}

Let us first prove a technical lemma, using the following well-known fact.

\lem[sect]{Let $L$ be a line bundle on a complex curve $X$, with $f\,:\,\O\to L$ be a non-zero section. Suppose it has zeros $p_1,\ldots,p_k$ with multiplicities $n_1,\ldots,n_k$. Then $c_1(L)\in\mathrm{H}^2_{\mathrm{DR}}(X)$ is Poincaré dual to the divisor $n_1p_1+\ldots,n_kp_k$.}

Let us call this divisor the \emph{zero-divisor} of $f$.

\defi{Let $B$ be a bundle of rank $r$ on a complex curve $X$. The {\it slope} of $B$ is:
$$\sl(B)=\frac{1}{r}\intg{X}{}c_1(B).$$}

\cor[hom]{Let $L,M$ be line bundles on a complex curve $X$. If $\sl(L)>\sl(M)$, then $\Hom(L,M)=0$. If $\sl(L)=\sl(M)$, then any $f\in\Hom(L,M)$ is either zero or an isomorphism.}

\demo{Let $f\in\Hom(L,M)$ neither zero, nor an isomorphism. It induces a non-zero section $\tilde{f}\,:\,\O\to M \otimes L^{-1}$. By Lemma \ref{sect}, if ${\cal D}$ is the zero-divisor of $\tilde{f}$,
$$\sl(M \otimes L^{-1})=\deg({{\cal D}})\le 0,$$
so $\sl(M)\le\sl(L)$.}

\lem[unique]{Let $X$ be a complex curve, $B$ a bundle of rank 2 on $X$ with splittings $B=L_1\oplus L_2=M_1\oplus M_2$. Then there are $i,j\in\{1,2\}$ such that $L_i\simeq M_j$ and $L_j\simeq M_i$. 
If $L_1\not\simeq L_2$, then one of these two isomorphisms is an equality, and if moreover $c_1(L_1)=c_1(L_2)$, both isomorphisms are equalities.}

\demo{Without loss of generality, $\sl(M_1)\ge\sl(L_1)\ge\sl(L_2)\ge\sl(M_2)$. 

If $\Hom(M_1,L_1)=\Hom(M_1,L_2)=0$, then $\Hom(M_1,M_1\oplus M_2)=\Hom(M_1,B)=0$, contradiction! Hence, by Corollary \ref{hom}, either $\sl(M_1)=\sl(L_1)$ and $M_1\simeq L_1$, or $\sl(M_1)=\sl(L_1)=\sl(L_2)$ and $M_1\simeq L_2$. In this second case, let us reverse the names of $L_1,L_2$, so that $M_1,L_1$ are anyway isomorphic.

We can not have $\Hom(L_1,M_2)=\Hom(L_2,M_2)=0$ either. Hence, either $L_2$ and $M_2$ are isomorphic (and we are done), or $\Hom(L_2,M_2)=0$ and $L_1,M_1,M_2$ are isomorphic. In this last case, $\Hom(L_2,M_1)=\Hom(L_2,M_2)=0$, so $\Hom(L_2,B)=0$, contradiction!

Now let us prove the second statement of Lemma \ref{unique}. We assume $\sl(M_1)=\sl(L_1)\ge\sl(L_2)=\sl(M_2)$ and $M_1\simeq L_1,M_2\simeq L_2$. If $\Hom(M_1,L_2)=0$, then the morphism $$f\,:\,M_1\hookrightarrow M_1\oplus M_2=L_1\oplus L_2\overset{p}{\to} L_2$$ induced by the projection $p$ to $L_2$ parallel to $L_1$ is zero, so $M_1=L_1$ and we are done. Else, $\sl(M_1)=\sl(L_2)$ and $M_1,L_2$ are isomorphic, so $L_1,L_2$ are isomorphic.

The third statement of Lemma \ref{unique} is proven similarly.}

Now, we can prove what we want.

\prop{Consider $X$ a $PK^0$ surface which is a ramified cover of a torus over three lines $C_1,C_2,C_3$ with no common intersection point and of distinct complex directions, with respective ramification orders $n_1\ge n_2\ge n_3$: $\pi\,:\, X\to \T$. Then the cotangent bundle $\Omega^1(X)$ does not split.}

\demo{Assume by contradiction $\Omega^1(X)=M_1\oplus M_2$. Let the line bundles $L_1,L_2,L_3$ be the complex directions of $\pi^{-1}(C_1),\pi^{-1}(C_2),\pi^{-1}(C_3)$ (defined on $\Omega^1(X\setminus X_s)$).

Let $C'_1\subset X$ be a parallel to $\pi^{-1}(C_1)$. By Lemmata \ref{localsplit} and \ref{hartogs}, $L_2,L_3$ can be extended to be defined globally on an open neighbourhood of $C'_1$, so that $\Omega^1(X)|_{C'_1}=L_2|_{C'_1}\oplus L_3|_{C'_1}$. Moreover, since $C'_1$ has positive intersection with $\pi^{-1}(C_2)$ and $\pi^{-1}(C_3)$ (and a generic choice of $C'_1$ gives generic intersection points), $L_2|_{C'_1}$ (of degree $n_2$) and $L_3|_{C'_1}$ (of degree $n_3$) are not isomorphic. According to Lemma \ref{unique}, $L_2|_{C'_1}=M_i|_{C'_1}$ for some $i$, say $i=2$. Passing to any parallel $C'_1$, we get $$L_2|_{X\setminus \pi^{-1}(C_1)}=M_2|_{X\setminus \pi^{-1}(C_1)}.$$

Similarly, $L_1|_{X\setminus \pi^{-1}(C_2)}=M_j|_{X\setminus \pi^{-1}(C_2)}$ for some $j$, and $j=1$ since $L_1,L_2$ are locally linearly independent as subbundles of $\Omega^1(X\setminus X_s)$.

Hence, the special tensor $\omega$ given by the splitting $\Omega^1(X)=M_1\oplus M_2$ satisfies: $$\omega_p=\alpha(p)\frac{\de z_1\otimes\de z_2+\de z_2\otimes\de z_1}{\de z_1\wedge\de z_2}$$ for all $p\in z\setminus\pi^{-1}(C_1\cup C_2)$ and some holomorphic never-vanishing function $\alpha$ on $X\setminus\pi^{-1}(C_1\cup C_2)$.

Let us write $\omega$ in local coordinates in a neighbourhood $U_p$ of a point $p\in\pi^{-1}(C_3)\setminus\pi^{-1}(C_1\cup C_2)$. Let $z_1,z$ be local coordinates in $U_z$ such that $L_1(U_p)=\C\cdot\de z_1,\, L_3(U_p)=\C\cdot\de z$, there is a chart given in these coordinates by: $$\phi_p\,:\, (z_1,z)\in U_p\mapsto (z_1,z^n)\in\C^2,$$ yielding $\de z_1 = \de z_1, nz^{n-1}\de z =\de z_3$. Since $\de z_3$ spans the complez direction of $C_3$ in $\T$, $\de z_2=a\de z_1+b\de z_3$ for some $a,b\in\C^*$. So, for all flat chart neighbourhood $V\subset z$, it holds in $U_p\cap V$ that $$\omega_{q=(z_1,z)}=\frac{a\alpha(q)}{nz^{n-1}}\frac{\de z_1\otimes\de z_1}{\de z_1\wedge\de z}+b\alpha(q)\frac{\de z_1\otimes\de z}{\de z_1\wedge\de z},$$
which is not bounded when $z$ tends to zero, whereas $\omega$ should also be defined at $p=(0,0)$, contradiction.}

\section*{Acknowledgements}

We are very grateful to Misha Verbitsky for proposing us this nice question, for explaining us so many mathematics on and around the subject of $PK$ manifolds and complex manifolds more generally, giving references as well as patiently exposing arguments that were not to be found in the literature.

\bibliographystyle{plain}
\bibliography{BibliographyV1} 

\end{document}